\documentclass[11pt]{article}
\usepackage[T1]{fontenc}
\usepackage{amsmath,amssymb,amsfonts}        
\usepackage{mathptmx}
\usepackage[scaled=0.91]{helvet}
\usepackage{courier}
\usepackage[dvips]{graphicx}
\usepackage{psfrag}
\usepackage{color}
\usepackage[numbers]{natbib}
\usepackage{setspace}

\title{A Finite Difference Method on
Quasi-uniform Mesh for Time-Fractional Advection-Diffusion Equations with Source Term}
\author{Riccardo Fazio and Alessandra Jannelli \\
\ \\
{\footnotesize Department of Mathematical and Computer Sciences,}\\
{\footnotesize Physical Sciences and Earth Sciences, University of Messina}\\
{\footnotesize Viale F. Stagno d'Alcontres 31, 98166 Messina, Italy}\\
{\footnotesize Email: rfazio@unime.it, ajannelli@unime.it} }
\date{}

\begin{document}

\maketitle

\begin{abstract}
The present paper deals with the numerical solution of time-fractional advection-diffusion 
equations involving the Caputo derivative with source term by means of an unconditionally stable implicit finite 
difference method on quasi-uniform grids. We use a special quasi-uniform mesh in order to improve the numerical 
accuracy of the classical discrete fractional formula for the Caputo derivative. 
The stability and the convergence of the method are discussed.
The error estimates established for a quasi-uniform grid and a uniform one are reported to support the theoretical results.
Numerical experiments are carried out to demonstrate the effectiveness of the method.
\end{abstract}


\section{Introduction}
The fractional partial differential equations (FPDEs) have become increasingly popular in recent years.
The interest in these equations comes from their mathematical structure and from their applications.
FPDEs have applied in various areas of engineering, science, finance, applied mathematics, bio-engineering and so on.
Though several ways to solve them theoretically have been proposed \cite{Oldham:TFC:1974,Kibas,Podlubny99},
including Green function method and Laplace and Fourier transform method  \cite{Samko,Miller},
the Adomian decomposition \cite{Daftardar-Gejji1,Cheng}, the Homotopy Perturbation Methods \cite{He2,Momani},
generally, numerical solution techniques are preferred when dealing with fractional models since the analytical solutions
are available for a few simple cases. In recent years, efficient numerical methods have been developed 
to solve fractional differential equation, including finite difference methods 
\cite{Meerschaert:FDA:2004,Tadjeran:ASA:2006,Liu:NST:2013,Ren:CDS:2013,Chen:INM:2013}, 
the finite volume method \cite{Hejazi:SCF:2014}, the finite element method \cite{Ervin:NAT:2007,Zhang:FDE:2012,Zeng:UFE:2013}  
and the spectral method \cite{Lin:FDS:2007,Doha:SJS:2013}.

As also for the non-fractional differential equations, finite difference methods
are one of the most important classes of numerical methods for solving FPDEs.
Zhuang and Liu \cite{Zhuang:IDA:2006} obtained an implicit difference approximation to solve 
the time-fractional diffusion equations.
Lin and Xu \cite{Lin:FDS:2007,Lin:FDS:2011} proposed the numerical solution by finite/difference 
approximations for a time-fractional diffusion equation. 
Liu et al. \cite{Liu:SCD:2007} developed an explicit difference method and an implicit difference method for solving a space-time fractional advection dispersion equation on a finite domain.
In \cite{Li:HOA:2009,Cao:HOA:2015,Cao:HOA:2016}, high order numerical difference schemes were constructed 
in order to solve the Caputo-type advection-diffusion equations.
Zhang et al. \cite{Zhang:FDM:2014} obtained a finite difference method for FPDEs involving the Caputo derivative 
on a non-uniform mesh.
Recently Jannelli et al. \cite{Jannelli:ANS:2017,Jannelli:ENS:2017} determined exact and numerical solutions
for the time-fractional advection-diffusion differential equations involving Riemann-Liouville derivative
with a non linear source term by means the Lie symmetries. They transformed the fractional partial differential 
equation into a fractional ordinary differential equations, which is then solved using the implicit backward differentiation formulas. 

The main goal of this paper is to construct an unconditionally stable implicit finite difference method 
for solving the time-fractional advection-diffusion equations (FADEs) with a nonhomogeneous source term  
involving the Caputo fractional derivative on quasi-uniform grids.
We choose to use a special quasi-uniform mesh in order to improve the numerical 
accuracy of the classical discrete fractional formula for the Caputo derivative, since the fractional derivatives are integrals with 
weakly singular kernel and the discretization on the uniform mesh may lead to poor accuracy. 
The consistency, stability and convergence of the proposed difference method are investigated. 
Three numerical examples are given to show the reliability and efficiency of the derived difference method.

\section{The mathematical model}
We consider the following linear time-fractional advection-diffusion equation
\begin{eqnarray}
 \frac{\partial^\alpha}{\partial t^\alpha} u(x,t) + K_1 \frac{\partial}{\partial x} u(x,t) - K_2 \frac{\partial^2}{\partial x^2} u(x,t) = f(x,t),  \quad  a < x < b , \quad 0 < t \leq T \ ,   \label{eq_0}
\end{eqnarray}
with the initial and boundary conditions given by
\begin{eqnarray*}
&& u(x,0) = \phi(x) \ ,  \qquad a \leq x \leq b \ , \\
&& u(a,t) = \varphi(t) \ , \qquad u(b,t) = \psi(t)\ , \qquad \qquad  0 < t \leq T \ , \nonumber
\end{eqnarray*}
where $u$ is the field variable that can represent, for example, the solute concentration, 
$K_1$ and  $K_2$ are the constant fluid velocity and the dispersion
coefficient, respectively. The time fractional derivative ${\displaystyle \frac{\partial^\alpha}{\partial t^\alpha} u(x,t)}$ is the $\alpha$ order Caputo fractional derivative defined by
\begin{equation}\label{op_c}
\frac{\partial^\alpha}{\partial t^\alpha} u(x,t) =\frac{1}{\Gamma (1-\alpha)} \int^t_0 \frac{\partial}{\partial s} u(x,s) (t-s)^{-\alpha}d s
 \qquad 0 < \alpha < 1  \ .
\end{equation}
The function $f(x,t)$ can be used to represent sources and sinks. $\phi(x)$, $\varphi(t)$ and $\psi(t)$ are known smooth functions.
We take $K_1 \neq 0$ and  $K_2>0$ and we assume the problem (\ref{eq_0}) has a unique and sufficiently smooth solution under the above
initial and boundary conditions.

The fractional equation (\ref{eq_0}) has been treated by a number of authors. It is presented as a useful approach
for the description of transport dynamics in complex systems which are governed by anomalous diffusion and
non-exponential relaxation patterns \cite{Metzler:RWG:2000}. The FADE is also used in groundwater hydrology research to
model the transport of passive tracers carried by fluid flow in a porous medium \cite{Benson:AFA:2000}.

It is to note that, when $\alpha = 1$, the model (\ref{eq_0}) reduces to the classical advection-diffusion equation with source term
\begin{eqnarray*}
\frac{\partial}{\partial t} u(x,t) + K_1 \frac{\partial}{\partial x} u(x,t) - K_2 \frac{\partial^2}{\partial x^2} u(x,t) = f(x,t), 
\qquad  a \leq x \leq b , \qquad 0 < t \leq T \ ,
\end{eqnarray*}
used in order to describe several phenomena of relevant interest in many fields of applied sciences.
In \cite{Jannelli:2003:MMP} and \cite{Jannelli:2009:SON}, a fractional step approach with variable time step is used in order to solve numerically mathematical models that describe evolution problems on a three-dimensional domains

\section{Discretization in time on a quasi-uniform mesh}
The main goal of this work is to construct an unconditionally stable implicit finite difference method defined on a quasi-uniform mesh.
In general, the existence of a weakly singular kernel $(t-s)^{-\alpha}$, $0 <\alpha<1$, in fractional derivatives makes it more difficult to get a higher-order scheme. Particularly when the solutions are not suitably smooth, numerical methods on uniform meshes seem to have a poor convergent rate. For these reasons, numerical schemes on non-uniform meshes have been developed in the last years.

In this Section, first, we construct the quasi-uniform mesh and then, in order to approximate the Caputo derivative on quasi-uniform grid, we define a suitable discrete fractional derivative formula.

We divide the interval $[0,T]$ into $N$ subintervals $[t^{n-1},t^{n}]$, for $n=1,\cdots,N$ and with $0 = t^0 < t^1 < \cdots t^N = T$.
We denote the time step $\Delta t^n $ as 
\begin{eqnarray*}
\Delta t^n= t^n - t^{n-1}, \qquad 1 \leq n \leq N \ ,
\end{eqnarray*}
and let 
\begin{eqnarray*}
\Delta t^{max} = \max_{1 \leq n \leq N} \Delta t^n \qquad \Delta t^{min} = \min_{1 \leq n \leq N} \Delta t^n \ .
\end{eqnarray*}

A sequence of mesh is called quasi-uniform if there exists a finite constant $\beta$ such that
\begin{eqnarray*}
\Delta t^{max} / \Delta t^{min} \leq \beta \ .
\end{eqnarray*}
In this case, it holds that $\Delta t^{max} \leq \beta T N^{-1}$. When $\beta = 1$, it hods that $\Delta t^n = T N^{-1}$
for all $n=1,\cdots,N$ and the mesh is reduced to a uniform mesh. 

A sequence of meshes is not quasi-uniform if 
\begin{eqnarray*}
\Delta t^{max} / \Delta t^{min} \rightarrow +\infty  \qquad \mbox{as} \qquad   N \rightarrow +\infty \ .
\end{eqnarray*}
The non-uniform mesh and quasi-uniform mesh methods have been used for solving differential equations by several authors.
In \cite{Fazio:FDS:2014,Fazio:BVP:2017}, quasi-uniform meshes were used for solving a class of boundary values problems on infinite domains.

\begin{figure}[!ht]
\centerline{
\psfrag{t}{$t$}
\includegraphics[width=.65\textwidth]{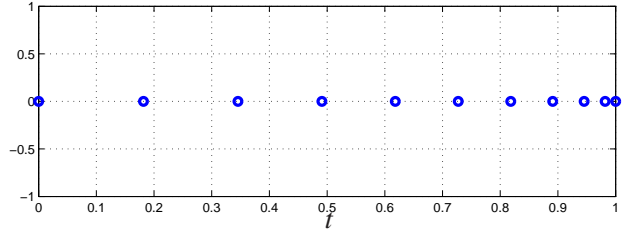}}
\caption{Quasi-uniform grid for $N=10$}
\label{grid_NU}
\end{figure}

In this work, we are interested on the quasi-uniform mesh \cite{Zhang:FDM:2014} defined as follows
\begin{eqnarray}\label{dt}
\Delta t^n =(N+1-n) \mu, \qquad  1 \leq n \leq N \ ,
\end{eqnarray}
where $\mu = {\displaystyle \frac{2T}{N(N+1)}}$. It is to note that the time steps $\{\Delta t^n\}_{n=1}^N$ are a monotonically decreasing sequence 
with $\Delta t^1 = O(N^{-1})$  and $\Delta t^N = O(N^{-2})$. Figure (\ref{grid_NU}) shows a sample of the quasi-uniform grid  (\ref{dt}) 
obtained for $N=10$.

To motivate the choice of the quasi-uniform mesh (\ref{dt}), we discretize the Caputo derivative (\ref{op_c}) of a function
$v(t)$ by means of the following classical approximate formula
\begin{eqnarray}\label{form_c}
\frac{d^\alpha}{d t^\alpha}  v(t^n) &=& \frac{1}{\Gamma (1-\alpha)} \int^{t^n}_0 v'(s) (t^n-s)^{-\alpha}d s  \\
&=&\frac{1}{\Gamma (1-\alpha)}  \sum_{k=1}^{n} \frac{v(t^k) - v(t^{k-1}) }{t^{k}-t^{k-1}} \int_{t^{k-1}}^{t^k} (t^n-s)^{-\alpha} d s 
+ R^n \nonumber
\end{eqnarray}
that is the well-known so-called L1 formula defined in \cite{Oldham:TFC:1974} where $R^n$ is the local truncation error. 
The L1 formula has been used for solving the fractional differential equations with Caputo derivatives 
(see \cite{Zhuang:IDA:2006,Ga0:CDS:2011}). Moreover, using the relationship between Caputo derivative 
and Riemann-Liouville fractional derivative, the L1 formula was also applied to time
fractional diffusion equation with Riemann-Liouville fractional derivative (see \cite{Langlands:ASA:2005,Zhang:EEC:2011}). High-order approximations such as compact difference scheme \cite{Du:CDS:2010,Ga0:CDS:2011,Zhang:EEC:2011,Chen:NSH:2010}
and spectral method \cite{Lin:FDS:2007,Lin:FDS:2011,Li:STS:2009} were applied to improve the spatial accuracy of
fractional diffusion equations. 
It is important to note that it is rather difficult to get a high-order time approximation due to the singularity of
fractional derivatives.

For any temporal mesh, for $0 < \alpha < 1 $ and $v(t) \in C^2[0,T]$, it can be verified that \cite{Zhang:FDM:2014}
\begin{eqnarray*}
\int_{0}^{t^n} v'(s) (t^n-s)^{-\alpha} d s = \sum_{k=1}^{n} \frac{v(t^k)- v(t^{k-1})}{t^{k}-t^{k-1}} \int_{t^{k-1}}^{t^k} (t^n-s)^{-\alpha} d s + r^n
\end{eqnarray*}
where
\begin{eqnarray*}
|r^n| \leq  \left( \frac{(\Delta t^n)^2}{2(1 -\alpha)} + \frac{(\Delta t_{max})^2}{8}\right) 
(\Delta t^n)^{-\alpha} \max_{0 \leq t \leq t^n} |v''(t)| \qquad 1 \leq n \leq N \ .
\end{eqnarray*}
Since 
\begin{eqnarray*} 
R^n &=& \frac{1}{\Gamma (1-\alpha)} \left [ \int^{t^n}_0 v'(s) (t^n-s)^{-\alpha}d s  -
     \sum_{k=1}^{n} \frac{v(t^k) - v(t^{k-1}) }{t^{k}-t^{k-1}} \int_{t^{k-1}}^{t^k} (t^n-s)^{-\alpha} d s  \right ] \nonumber \\
		&=&	\frac{1}{\Gamma (1-\alpha)} r^n \ ,
\end{eqnarray*}
we obtain 
\begin{eqnarray*}
|R^n| \leq  \frac{1}{\Gamma (1-\alpha)} \left( \frac{(\Delta t^n)^2}{2(1 -\alpha)} + \frac{(\Delta t_{max})^2}{8}\right) 
(\Delta t^n)^{-\alpha} \max_{0 \leq t \leq t^n} |v''(t)| \qquad 1 \leq n \leq N \ .
\end{eqnarray*}
For the uniform mesh, {\it i.e.}, $\Delta t^n = \Delta t $ for all $n = 1, 2,\cdots,N$, then $R^n = O(\Delta t^{2 -\alpha})$ 
(see \cite{Sun:FDD:2006,Lin:FDS:2007}). 
From the truncation error estimate of the L1 formula, it is clear that the accuracy is dependent on
the fractional order $\alpha$. This is justified since a weakly singular kernel $(t - s)^{-\alpha} $ is contained in the integral. 

In order to improve the accuracy of the L1 numerical approximation of derivative of fractional order,
it is possible to consider non-uniform meshes. Numerical methods developed with the non-uniform meshes have been developed 
for solving integro-differential equations with weakly singular. Mustapha \cite{Mustapha:IFD:2011}, Yuste and Quintana-Murillo 
\cite{Yuste:FDS:2012,Yuste:FAR:2016} proposed an implicit finite-difference time-stepping method for discretizing of the time diffusion equation. 

Zhang et al. \cite{Zhang:FDM:2014} obtained a numerical integration formula for any $ \alpha \in (0, 1)$ by employing the special quasi--uniform grid (\ref{dt}) and obtained the following results: 
for the quasi-uniform mesh (\ref{dt}), for $0 < \alpha < 1 $ and $v(t) \in C^2[0,T]$, it holds that 
\begin{eqnarray*}
\int_{0}^{t^n} v'(s) (t^n-s)^{-\alpha} d s = \sum_{k=1}^{n} \frac{v(t^k)- v(t^{k-1})}{t^{k}-t^{k-1}} \int_{t^{k-1}}^{t^k} (t^n-s)^{-\alpha} d s + r^n \quad 1 \leq n \leq N
\end{eqnarray*}
where
\begin{eqnarray*}
&&| r^n | \leq  \left (1 + \alpha + \frac{2^{1-\alpha}}{1-\alpha} \right)  \max_{0 \leq t \leq t^n} |v''(t)|  T^{2-\alpha} (N+1)^{\alpha-2} \ ,
\quad 1 \leq n \leq N-1 \ ,
\end{eqnarray*}
and 
\begin{eqnarray*}
\int_{0}^{t^N} v'(s) (t^N-s)^{-\alpha} d s = \sum_{k=1}^{N} \frac{v(t^k)- v(t^{k-1})}{t^{k}-t^{k-1}} \int_{t^{k-1}}^{t^k} (t^N-s)^{-\alpha} d s + r^N 
\end{eqnarray*}
where
\begin{eqnarray*}
&&| r^N | \leq  \frac{1+\alpha}{1-\alpha} 2^{1-\alpha}  \max_{0 \leq t \leq T} |v''(t)|  T^{2-\alpha} N^{-2} \ . 
\end{eqnarray*}
Then, we obtain
\begin{eqnarray}\label{Error}
&&| R^n | \leq \frac{1}{\Gamma (1-\alpha)} \left (1 + \alpha + \frac{2^{1-\alpha}}{1-\alpha} \right)  \max_{0 \leq t \leq t^n} |v''(t)|  T^{2-\alpha} (N+1)^{\alpha-2} \ ,
\quad 1 \leq n \leq N-1 \ , \nonumber \\ 
\\
&&| R^N | \leq \frac{1}{\Gamma (1-\alpha)} \frac{1+\alpha}{1-\alpha} 2^{1-\alpha}  \max_{0 \leq t \leq T} |v''(t)|  T^{2-\alpha} N^{-2} \ . \nonumber
\end{eqnarray}
See \cite{Zhang:FDM:2014} for more details. We use these estimates to study the consistency, stability and 
convergence of the method presented in the following sections.

\section{An implicit finite difference method}
For the derivation of the implicit difference method, first, we construct a computational uniform grid in the $x$ direction, that is the spatial size of the mesh $\Delta x = x_j - x_{j-1}$ is constant, for $1 \leq j \leq J$, and quasi-uniform in the time direction.
We define the mesh points $(x_j,t^n)$ with $x_j= a + j \Delta x$, $j=0,\cdots,J$ and $t^n = t^{n-1} + \Delta t^n$, for $n=1,\cdots,N$,
with $\Delta t^n$ defined by (\ref{dt}). $J$ and $N$ are positive integers. 
We denote by $U_j^n$ the numerical approximation provided by the difference method of the exact solution $u(x_j,t^n)$ 
at the mesh points $(x_j,t^n)$, for $j=0,\cdots,J$ and $n=0,\cdots,N$.

As usual, we discretize the first $\partial / \partial x$ and second order $\partial^2/\partial x^2$ spatial derivatives by means of the 
second order three-point central difference formula so that 
\begin{eqnarray}
&& \frac{\partial}{\partial x} u (x_j,t^n) =  \frac{u(x_{j+1},t^n) - u(x_{j-1},t^n) }{2 \Delta x} + O(\Delta x^2) \label{d1} \\
\nonumber \\
&& \frac{\partial^2}{\partial x^2} u (x_j,t^n) =  \frac{u(x_{j+1},t^n) - 2 u(x_{j},t^n) + u(x_{j-1},t^n) }{\Delta x^2} + O(\Delta x^2) 
\label{d2} \ .
\end{eqnarray}
According to the discretization for the Caputo derivative (\ref{form_c}), we approximate the time fractional derivative in 
(\ref{eq_0}) as follows
\begin{eqnarray}
&& \frac{\partial^\alpha}{\partial t^\alpha}  u(x_j,t^n) =  \label{Error1} \\
&& \frac{1}{\Gamma (2-\alpha)}  \sum_{k=1}^{n} \frac{u(x_j,t^k) - u(x_j,t^{k-1})}{t^{k}-t^{k-1}} 
\left[ (t^n-t^{k-1})^{1-\alpha} - (t^n-t^{k})^{1-\alpha} \right]  + R^n_j  \nonumber \ ,
\end{eqnarray}
where 
\begin{eqnarray*}
\frac{1}{\Gamma(1-\alpha)} \int_{t^{k-1}}^{t^k} (t^n-s)^{-\alpha} ds = \frac{1}{\Gamma(2-\alpha)} [(t^n-t^{k-1})^{1-\alpha} - (t^n-t^{k})^{1-\alpha}] \ .
\end{eqnarray*} 
Replacing $u(x_j,t^n)$ with its numerical approximation $U^n_j$ and neglecting the local truncation errors, the time-fractional advection-diffusion equation (\ref{eq_0}) is discretized as follows
\begin{eqnarray}\label{meth_0}
\frac{1}{\Gamma (2-\alpha)} \sum_{k=1}^{n} T_{n,k} (U_j^k - U_j^{k-1}) + K_1 \frac{U_{j+1}^n - U_{j-1}^n }{2 \Delta x} - K_2 \frac{U_{j+1}^n - 2 U_{j}^n +U_{j-1}^n }{\Delta x^2} = f^n_j 
\end{eqnarray}
for $ 1 \le n \le N $ and $1 \le j \le J-1$, where $f^n_j=f(x_j,t^n)$ is the source term and 
\begin{eqnarray*}
 T_{n,k} = \frac{(t^n-t^{k-1})^{1-\alpha} - (t^n-t^{k})^{1-\alpha} }{t^{k}-t^{k-1}} \qquad 1 \leq k \leq n \ , \qquad 1 \leq n \leq N  \ .
\end{eqnarray*} 
The initial and boundary conditions can be rewritten as
\begin{eqnarray}\label{Ic_Bc}
&& U^0_j = \phi(x_j)  \qquad 0 \leq j \leq J  \nonumber \\
&& U^n_0 = \varphi(t^n) \qquad U^n_J = \psi(t^n) \qquad 0 \leq n \leq N  \ .
\end{eqnarray}
For any temporal meshes on $[0,T]$ and for any $n=1,2,\cdots,N$, it holds that \cite{Zhuang:IDA:2006}
\begin{eqnarray*}
&& T_{n,k} > 0 \qquad \qquad T_{n,k} > T_{n,k-1} \ , \qquad 1 \leq k \leq n \ .
\end{eqnarray*}
Taking into account that $T_{n,n} = (\Delta t^n)^{-\alpha}$, we can write
\begin{eqnarray*}
&& \sum_{k=1}^{n} T_{n,k} (U_j^k - U_j^{k-1}) =\sum_{k=1}^{n-1} T_{n,k} (U_j^k - U_j^{k-1}) + (\Delta t^n)^{-\alpha}  (U_j^n - U_j^{n-1})
\end{eqnarray*}
then, we obtain the following implicit finite difference scheme
\begin{eqnarray}
&& (-K1-K2) \ U_{j-1}^n + (1 + 2 \ K2) \ U_j^n + (K1 - K2) \ U_{j+1}^n   \label{meth} \\
&& = U_j^{n-1} - (\Delta t^n)^{\alpha} \sum_{k=1}^{n-1} T_{n,k} (U_j^k - U_j^{k-1}) + F^n_j \ ,  \qquad 1 \le n \le N, \qquad 1 \le j \le J-1 \ ,
\nonumber 
\end{eqnarray}
where we set
\begin{eqnarray*}
K1= \frac{K_1 (\Delta t^n)^{\alpha}  \Gamma(2-\alpha)}{2 \Delta x} \ ,  \quad K2 = \frac{K_2 (\Delta t^n)^{\alpha} \Gamma (2-\alpha)}{\Delta x^2} \ , \quad F_j^n = (\Delta t^n)^{\alpha} \Gamma(2-\alpha) f^n_j  \ .
\end{eqnarray*}
The equation (\ref{meth}) can be written in vectorial form as
\begin{eqnarray}
&& K U^n = \mathcal{L} U^{n-1} + F^n  \label{meth_m} \ ,
\end{eqnarray}
where $K$ is a tridiagonal matrix
and where we denote with $ \mathcal{L}$ the following difference operator
\begin{eqnarray}\label{L}
&& \mathcal{L} \ U^{n-1} =  U_j^{n-1} - (\Delta t^n)^{\alpha} \sum_{k=1}^{n-1} T_{n,k} (U_j^k - U_j^{k-1})  \ .
\end{eqnarray}
Here and in the following we assume the convention that the summation is equal to zero if the lower bound is larger that the upper bound. 
The obtained method (\ref{meth_m}) is implicit, in order to compute the numerical solution $U^n$ 
a system with the tridiagonal coefficients matrix $K$ has to be solved.

It is interesting to note that the operator $\mathcal{L}$ is a kind operator with memory, due to the non-local character of the fractional derivative, this means that the effect on $U$ at time $t^n$, $U^n$, depends on all the previous values, $ U^0,U^1,\cdots,U^{n-1}$, evaluated at all the previous time $t^0,t^1,\cdots,t^{n-1}$. 

The main difference compared to the non-fractional case is that, in order to evaluate $\mathcal{L}$, 
the numerical solutions for all the $n$ previous time values $t^0,t^1,\cdots,t^{n-1}$ are required, 
while for non-fractional equations, only the solution to the previous value $t^{n-1}$ is used.
The computational cost to compute the solution at the time $t^n$ from the solution at the time $t^{n-1}$
grows as $n$, that is, grows as the number of terms in the summation that compares in the second term of the (\ref{L}). 
This implies that the computational cost to go from $t^0$ to $t^n$ grows as $n^2$.


\section{Consistency, Stability and Convergence}
In this Section, we discuss the consistency, the stability and the convergence of the implicit finite difference scheme.

\bigskip

\noindent
{\bf Consistency}.
According to Eqs. (\ref{d1}), (\ref{d2}) and (\ref{Error1}), the local truncation error of the difference scheme (\ref{meth_0}) is 
\begin{eqnarray}
R^n_j &=& \frac{1}{\Gamma (2-\alpha)} \sum_{k=1}^{n} T_{n,k} (u(x_j,t^k) - u(x_j,t^{k-1}))  \nonumber \\
&& + K_1 \frac{u(x_{j+1},t^n) - u(x_{j-1},t^n) }{2 \Delta x} - K_2 \frac{u(x_{j+1},t^n) - 2 u(x_{j},t^n) + u(x_{j-1},t^n) }{\Delta x^2} - f(x_j,t^n)  \nonumber \\
&=& \frac{1}{\Gamma (2-\alpha)} \sum_{k=1}^{n} T_{n,k} (u(x_j,t^k) - u(x_j,t^{k-1})) -\frac{\partial^\alpha}{\partial t^\alpha}  u(x_j,t^n)  \nonumber \\
&& + K_1 \left [\frac{u(x_{j+1},t^n) - u(x_{j-1},t^n) }{2 \Delta x} - \frac{\partial}{\partial x} u (x_j,t^n) \right ]  \nonumber \\
&& - K_2 \left [\frac{u(x_{j+1},t^n) - 2 u(x_{j},t^n) + u(x_{j-1},t^n) }{\Delta x^2} - \frac{\partial^2}{\partial x^2} u (x_j,t^n)  \right ]   \nonumber \\
&=& O(N^{\alpha-2}) + K_1 O(\Delta x^2) + K_2 O(\Delta x^2) = O(N^{\alpha-2} + \Delta x^2) \ .
\end{eqnarray}
The implicit finite difference scheme defined by (\ref{meth_0}) or (\ref{meth}) is consistent with the model (\ref{eq_0}) of order 
$O(N^{\alpha-2} + \Delta x^2)$.

\bigskip

\noindent
{\bf Stability}. For the stability analysis of the implicit finite difference scheme
we rewrite the equations (\ref{meth}) as
\begin{eqnarray}\label{meth_stb} 
&& (-K1-K2) \ U_{j-1}^n + (1 + 2 \ K2)  U_j^n + (K1 - K2) \ U_{j+1}^n  \\
&& = (\Delta t^n)^{\alpha} \sum_{k=0}^{n-1} (T_{n,k+1} - T_{n,k}) U_j^k + F^n_j \ , \qquad 1 \le n \le N \ ,\qquad 1 \le j \le J-1  \ .\nonumber
\end{eqnarray}
\noindent
Let $\bar{U}^n_j$ be another approximate solution of the difference scheme (\ref{meth}), and let 
\begin{eqnarray*}
\rho^n_j =  U^n_j - \bar{U}^n_j \ , \qquad 0 \leq j \leq J \ , \qquad 1 \leq n \leq N \ ,
\end{eqnarray*}
be the corresponding round-off error. We let
\begin{eqnarray*}
\rho^n = (\rho_0^n, \rho_1^n, \cdots, \rho_{J}^n)^T \ ,
\end{eqnarray*}
and we consider the infinity norm\begin{eqnarray*}
|| \rho^n ||_\infty = \max_{0 \leq j \leq J}  | \rho^n_j| = | \rho^n_i| \ .
\end{eqnarray*}
The round-off error satisfies the following round-off equations 
\begin{eqnarray}\label{meth_rho}
&& (-K1-K2) \rho_{j-1}^n + (1 + 2 \ K2) \rho_j^n + (K1 - K2)\rho_{j+1}^n \\
&&  = (\Delta t^n)^{\alpha} \sum_{k=0}^{n-1} (T_{n,k+1} - T_{n,k}) \rho_j^k  \ . \nonumber
\end{eqnarray}
In order to check whether the finite difference scheme is stable, we study how the size of the round-off error $\rho^n$ evolves in time.

\noindent
We define
\begin{eqnarray*}
&& L_1 \rho^n_j = (-K1-K2) \rho_{j-1}^n + (1 + 2 \ K2) \rho_j^n + (K1 - K2)\rho_{j+1}^n 
\end{eqnarray*}
and 
\begin{eqnarray}\label{L2}
&& L_2 \rho^{n-1}_{j} = (\Delta t^n)^{\alpha} \sum_{k=0}^{n-1} (T_{n,k+1} - T_{n,k}) \rho_j^k  \ .
\end{eqnarray}

\noindent
Eq (\ref{meth_rho}) can be written as
\begin{eqnarray*}
&& L_1 \rho^n_j = L_2 \rho^{n-1}_j \ .
\end{eqnarray*}

\noindent
From (\ref{L2}) and taking into account that $T_{n,k+1} - T_{n,k} >0$, we have
\begin{eqnarray}
| L_2 \rho^{n-1}_j | = |(\Delta t^n)^{\alpha} \sum_{k=0}^{n-1} (T_{n,k+1} - T_{n,k}) \rho_j^k |  \leq |\rho_{j}^{n-1} | (\Delta t^n)^{\alpha}  \sum_{k=0}^{n-1} (T_{n,k+1} - T_{n,k})  \label{l_max}
\end{eqnarray}
where we have defined 
\begin{eqnarray*}
|\rho_{j}^{n-1}  | = \max_{0 \leq k \leq n-1}  |\rho^k_j| \ . 
\end{eqnarray*}
But, since $ \sum_{k=0}^{n-1} (T_{n,k+1} - T_{n,k}) = (\Delta t^n)^{-\alpha} $ and
recalling that $T_{n,n}= (\Delta t^n)^{-\alpha} $ and $T_{n,0}=0$, it has
\begin{eqnarray*}
| L_2 \rho^{n-1}_j | = |(\Delta t^n)^{\alpha} \sum_{k=0}^{n-1} (T_{n,k+1} - T_{n,k}) \rho_j^k |  \leq |\rho_{j}^{n-1} |  \ ,
\end{eqnarray*}
Thus, we can conclude that 
\begin{eqnarray*}
|| \rho^n ||_\infty &=& | \rho^n_i|\  = \ | (-K1-K2) \rho_{i}^n + (1 + 2 \ K2) \rho_i^n + (K1 - K2)\rho_{i}^n  | \\
&=& | L_1 \rho^{n}_i | \ = \ | L_2 \rho^{n-1}_i | \ \leq \ |\rho_{i}^{n-1} | \ = \ || \rho^{n-1} ||_\infty  \ ,
\end{eqnarray*}
i.e. 
$$|| \rho^n ||_\infty \leq || \rho^0 ||_\infty \ ,  $$
for $1 \leq n \leq N$. This means that the present method is unconditionally stable.

\bigskip

\noindent
{\bf Convergence}.
Let $u(x_j, t_n)$ be the exact solution of Eqs.  (\ref{eq_0}) at mesh point $(x_j, t^n)$ for $j = 0,1,\cdots,J$ and $n = 0,1,\cdots,N$.
Denoting $\epsilon^n_j = u(x_j, t_n) -U_j^n$, we get the error equations
\begin{eqnarray}
&& L_1 \epsilon^n_j = L_2 \epsilon^{n-1}_j + R^n_j \ , \qquad j = 0,1,\cdots,J\ , \qquad  n = 0,1,\cdots,N  \ ,
\end{eqnarray}
with $\epsilon^0_j = 0$, for $j=0,1,\cdots,J$ and $\epsilon^n_J=0$, for $n=1,2,\cdots,N$. $R^n_j $ is the local truncation error.
We introduce the following norm
\begin{eqnarray*}
|| \epsilon^n ||_\infty = \max_{0 \leq j \leq J}  | \epsilon^n_j| = | \epsilon^n_i| \ ,
\end{eqnarray*}
then, we obtain
\begin{eqnarray*}
 || \epsilon^n ||_\infty &=& | \epsilon^n_i| \ = \ | L_1 \epsilon^{n}_i | \ = \ | L_2 \epsilon^{n-1}_i + R^{n}_i|   \\
&\leq& | L_2 \epsilon^{n-1}_i| \ + \ |R^{n}_i| \ \leq \ |\epsilon_{i}^{n-1} | \ + 
                                   \ |R^{n}_i| \ \leq \ || \epsilon^{n-1} ||_\infty + \ R_{max} \ ,
\end{eqnarray*}
where $R_{max} = \max_{n,i}  |R^n_i| $. 
For $n=1$, we have
\begin{eqnarray*}
|| \epsilon^1 ||_\infty \leq || \epsilon^0 ||_\infty + R_{max} = R_{max}\ .
\end{eqnarray*}
Then 
\begin{eqnarray*}
|| \epsilon^n ||_\infty \leq R_{max} \ , \qquad 1 \leq n \leq N \ .
\end{eqnarray*}
To obtain the error estimates of the numerical solutions, we need the uniform error bounds on all time levels. 
Then, applying the first of the estimates (\ref{Error}), we can
consider
\begin{eqnarray}
&& R_{max} \leq C_R (N^{\alpha-2} + \Delta x^2) \ , \qquad 1\leq n \leq N \ ,
\end{eqnarray}
where $C_R$ is a positive constant that is dependent on $T$, $\alpha$ and the exact solution $u(x,t)$, but independent of $N$ and $\Delta x$.
Then, we prove that the solution of the difference method (\ref{meth}), with initial and boundary conditions given by (\ref{Ic_Bc}),
is convergent.

\section{Numerical experiments}
In this Section, we report some numerical examples of the FPDEs to demonstrate the accuracy and efficiency of the numerical method.
The first examples are chosen in such a way the exact solution of FPDE can be evaluated analytically to show that
the approach proposed in this paper properly works. The exact solutions allow to verify the accuracy and the order of
convergence of the numerical solution.  We use the proposed method on the quasi-uniform grid and on a uniform grid and compare 
the numerical results, observing that the finite difference scheme generates more accurate numerical solutions on the quasi-uniform grid 
than on the uniform one. 
In the third test, we solve a FADE of physical interest with the source term chosen as a linear function of the field variable.
in order to show that the method is applicable to a wide class of FADEs.
By this last test, we illustrate how the changes in the solution behavior arise when the fractional order is varied.

\bigskip

\noindent
{\bf Example 1}. We consider the following FADE 
\begin{eqnarray}\label{test1}
&& \frac{\partial^\alpha}{\partial t^\alpha} u(x,t) + \frac{\partial}{\partial x} u(x,t) - \frac{\partial^2}{\partial x^2} u(x,t) 
= f(x,t),  \quad  0 < x < 1 , \quad 0 < t \leq T \ , \nonumber  \\
\nonumber  \\
&& u(x,0) = 0 \ ,   \qquad \qquad 0 \leq x \leq 1 \ , \\
&& u(0,t) = t^\beta \ ,\qquad \qquad u(1,t) = e t^\beta \ ,  \qquad  0 < t \leq T \ , \nonumber
\end{eqnarray}
where the source term is given by
\begin{eqnarray*}
&& f(x,t) = \frac{\Gamma(\beta+1)}{\Gamma(\beta+1-\alpha)}  e^x t^{\beta-\alpha} \ .
\end{eqnarray*}
The exact solution is found to be
\begin{eqnarray*}
&& u(x,t) =  e^x t^{\beta} \ .
\end{eqnarray*}

\begin{figure}[!hbt]
\centerline{
\psfrag{x}{{$x$}}
\psfrag{t}{{$t$}}
\psfrag{U_j^n,u(x_j,t^n)}{{$U_j^n, u(x_j,t^n)$}}
\includegraphics[width=.65\textwidth]{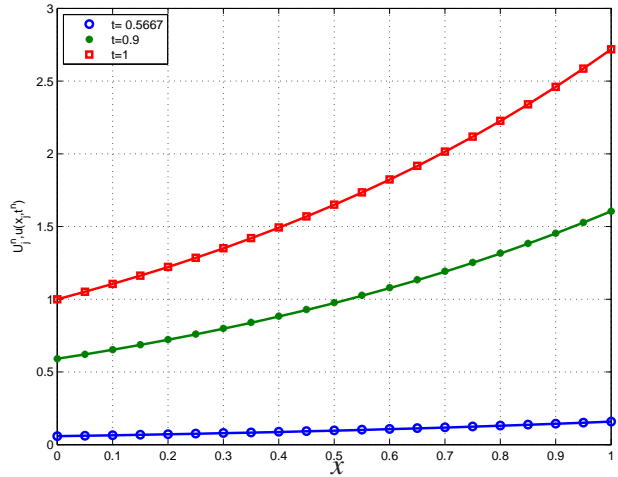}}
\caption{Comparison of the exact and the numerical solutions of the FPDE (\ref{test1}) for $\alpha=0.5$ with $N=J=20$ at different 
time steps. Solid lines: exact solution, star, circle and square points: numerical solution.}
\label{fig1}
\end{figure}

In this example, we take $K_1=K_2=1$  with $\alpha = 0.5$ and $\beta = 5$. 
Figure \ref{fig1} shows the comparison between the numerical solution $U^n_j$ and the exact solution $u(x_j,t^n)$ at different times computed with $N=J=20$. From the Figure \ref{fig1}, it can be seen that the numerical solution $U^n_j$ is in good agreement with the exact solution 
$u(x_j,t^n)$. The exact solution is reported with the solid line.

In order to investigate the temporal error and the convergence order of the numerical difference method,
we define the maximum error between the exact solution $u(x_j,t^N)$ and the numerical solution $U_j^N$ at the final time $t^N$
\begin{eqnarray}\label{error}
e_\infty(N,J) &=& \max_{1 \le j \le J} | u(x_j,t^N) - U_j^N|
\end{eqnarray}
and the convergence order as follows
\begin{eqnarray}\label{order}
Order= \log_2 \left( \frac{e_\infty(N,J)}{e_\infty(2N,J)} \right ) \ .
\end{eqnarray}

In order to show the efficiency of the method,  we define
${\bar U}_j^n$ the numerical solution obtained by the implicit difference method on a uniform grid defined by $\Delta t^n = T/N$, 
for $n=1,\cdots,N$. We use notations similar to the (\ref{error}) and (\ref{order}) to define the maximum error $\bar{e}_\infty(N,J)$  between the exact solution $u(x_j,t^N)$ and the numerical solution ${\bar U}_j^N$ at the final time $t^N$ and the convergence $Order_{\bar{U}}$.
In this test, we fix $J=100$, a value large enough such that the spatial error is negligible as compared with the temporal error.
Table~\ref{table1} shows the values of $e_\infty$ and $\bar{e}_\infty$ and the corresponding numerical convergence orders 
for $\alpha=0.1,0.5$ and $0.9$. It can be seen that the method is stable and convergent for solving the problem (\ref{test1}) on both
the computational grids. 
By using the L1 formula on quasi-uniform mesh, we improve the order of accuracy in time of the proposed method.
In fact, we observe that the solutions are more accurate on the quasi-uniform mesh and the convergence order 
on the quasi-uniform mesh is greater than the convergence order obtained with the uniform mesh. 
The numerical results agree well with the theoretical results.

\begin{table}[!h]
\centering
\renewcommand\arraystretch{1.2}
\begin{tabular}{lrrrrrr}
\hline
$ \alpha $ & $N$ & $e_\infty$ & $Order$ & ${\bar e}_\infty$ & $Order_{\bar U}$ \\
\hline
$0.1$  & $10$  & $ 3.6363e-04 $ & $ $        & $ 1.3741e-03$ & \\
       & $20$  & $ 9.1021e-05 $ & $ 1.9982 $ & $ 4.3575e-04$ & $ 1.6570 $ \\
       & $40$  & $ 2.2054e-05 $ & $ 2.0452 $ & $ 1.3222e-04$ & $ 1.7206 $ \\
       & $80$  & $ 4.4649e-06 $ & $ 2.3043 $ & $ 3.8355e-05$ & $ 1.7854 $ \\
\\
$0.5$  & $10$  & $ 5.5793e-03 $ & $ $        & $ 1.9875e-02$ & \\
       & $20$  & $ 1.7121e-03 $ & $ 1.7044 $ & $ 7.7106e-03$ & $ 1.3660 $ \\
       & $40$  & $ 5.4236e-04 $ & $ 1.6584 $ & $ 2.8895e-03$ & $ 1.4160 $ \\
       & $80$  & $ 1.7544e-04 $ & $ 1.6283 $ & $ 1.0602e-03$ & $ 1.4465 $ \\
\\
$0.9$  & $10$  & $ 4.6556e-02 $ & $ $        & $ 1.0115e-01$ & \\
       & $20$  & $ 2.1358e-02 $ & $ 1.1242 $ & $ 4.9512e-02$ & $ 1.0307 $ \\
       & $40$  & $ 9.8735e-03 $ & $ 1.1131 $ & $ 2.3691e-02$ & $ 1.0634 $ \\
       & $80$  & $ 4.5790e-03 $ & $ 1.1085 $ & $ 1.1201e-02$ & $ 1.0807 $ \\
\hline
\end{tabular}
\caption{$e_\infty$, ${\bar e}_\infty$ and convergence orders, related to the numerical solutions $U_j^n$ and ${\bar U}_j^n$ respectively, 
for different values of $N$ and $\alpha$, with $J=100$.}
\label{table1}
\end{table}

\bigskip

\noindent
{\bf Example 2}. We consider the following FADE 
\begin{eqnarray}\label{test1_1}
&& \frac{\partial^\alpha}{\partial t^\alpha} u(x,t) + \frac{\partial}{\partial x} u(x,t) - \frac{\partial^2}{\partial x^2} u(x,t) 
= f(x,t),  \quad  0 < x < 1 , \quad 0 < t \leq T \ , \nonumber  \\
\nonumber  \\
&& u(x,0) = 0 \ ,   \qquad \qquad 0 \leq x \leq 1 \ , \\
&& u(0,t) = 0 \ ,\qquad \qquad u(1,t) = t^3 \ ,  \qquad  0 < t \leq T \ , \nonumber
\end{eqnarray}
where the source term is given by
\begin{eqnarray*}
&& f(x,t) = \frac{6}{\Gamma(4-\alpha)}  x^2 t^{3-\alpha} + 2 t^3 (x-1) 
\end{eqnarray*}
and the exact solution is 
\begin{eqnarray*}
&& u(x,t) =  x^2 t^3 \ .
\end{eqnarray*}

\begin{figure}[!hbt]
\centerline{
\psfrag{x}{{$x$}}
\psfrag{t}{{$t$}}
\psfrag{U_j^n,u(x_j,t^n)}{{$U_j^n, u(x_j,t^n)$}}
\includegraphics[width=.65\textwidth]{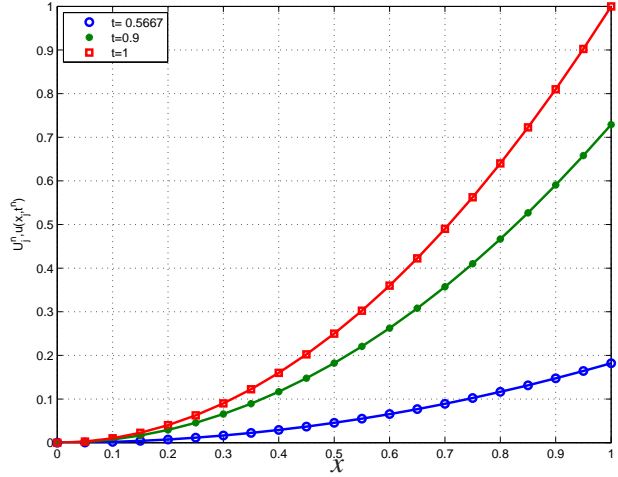}}
\caption{Comparison of the exact and the numerical solutions of the FPDE (\ref{test1_1}) for $\alpha=0.1$ with $N=J=20$ at different 
time steps. Solid lines: exact solution, star, circle and square points: numerical solution.}
\label{fig1_1}
\end{figure}

We take $K_1=K_2=1$ and set $\alpha = 0.1$. 
Figure \ref{fig1_1} shows the comparison between the numerical solution $U^n_j$ and the exact solution $u(x_j,t^n)$ at different times computed with $N=J=20$. The exact solution is reported with the solid line.
The numerical solution $U^n_j$ is in good agreement with the exact solution $u(x_j,t^n)$. 

The values of $e_\infty$ and $\bar{e}_\infty$ and the corresponding numerical convergence orders obtained
for $\alpha=0.1,0.5$ and $0.9$ are reported in the table~\ref{table1_1}. 
The results confirm that the method is stable and convergent for solving the problem (\ref{test1_1}) on both
the computational grids. The numerical solutions computed on quasi-uniform mesh are more accurate and the convergence order 
on the quasi-uniform mesh is greater than the convergence order obtained with the uniform mesh. 
The errors $e_\infty$ and $\bar{e}_\infty$ satisfy the relationships (\ref{error}).
The convergence orders $Order$ and  $Order_{\bar U}$ satisfy the relationships (\ref{order}).

\begin{table}[!h]
\centering
\renewcommand\arraystretch{1.2}
\begin{tabular}{lrrrrrr}
\hline
$ \alpha $ & $N$ & $e_\infty$ & $Order$ & ${\bar e}_\infty$ & $Order_{\bar U}$ \\
\hline
$0.1$  & $10$  & $ 3.7836e-05 $ & $ $        & $ 9.4723e-05$ & \\
       & $20$  & $ 9.6542e-06 $ & $ 1.9705 $ & $ 2.8745e-05$ & $ 1.7204 $ \\
       & $40$  & $ 2.4674e-06 $ & $ 1.9682 $ & $ 8.5307e-06$ & $ 1.7526 $ \\
       & $80$  & $ 6.2843e-07 $ & $ 1.9732 $ & $ 2.4911e-06$ & $ 1.7759 $ \\
\\
$0.5$  & $10$  & $ 4.4597e-04 $ & $ $        & $ 1.3182e-03$ & \\
       & $20$  & $ 1.3436e-04 $ & $ 1.7309 $ & $ 4.8998e-04$ & $ 1.4278 $ \\
       & $40$  & $ 4.1920e-05 $ & $ 1.6803 $ & $ 1.7898e-04$ & $ 1.4529 $ \\
       & $80$  & $ 1.3435e-05 $ & $ 1.6416 $ & $ 6.4671e-05$ & $ 1.4686 $ \\
\\
$0.9$  & $10$  & $ 3.3229e-03 $ & $ $        & $ 6.9875e-03$ & \\
       & $20$  & $ 1.4844e-03 $ & $ 1.1625 $ & $ 3.3290e-03$ & $ 1.0682 $ \\
       & $40$  & $ 6.7617e-04 $ & $ 1.1345 $ & $ 1.5524e-03$ & $ 1.0837 $ \\
       & $80$  & $ 3.1119e-04 $ & $ 1.1196 $ & $ 7.2847e-04$ & $ 1.0915 $ \\
\hline
\end{tabular}
\caption{$e_\infty$, ${\bar e}_\infty$ and convergence orders, related to the numerical solutions $U_j^n$ and ${\bar U}_j^n$ respectively, 
for different values of $N$ and $\alpha$, with $J=100$.}
\label{table1_1}
\end{table}

\bigskip

\noindent
{\bf Example 3}. 
In this test, in order to show the efficiency of the method,
we solve the following FADE 
\begin{eqnarray}\label{test2}
&& \frac{\partial^\alpha}{\partial t^\alpha} u(x,t) + K_1 \frac{\partial}{\partial x} u(x,t) - K_2 \frac{\partial^2}{\partial x^2} u(x,t) 
= f(u(x,t)) \ ,  \nonumber  \\
\nonumber  \\
&& u(x,0) = x^2(5-x)^2 \ , \qquad a \leq x \leq b \ , \\
&& u(a,t) = u(b,t) = 0 \ , \qquad  0 < t \leq T \ , \nonumber
\end{eqnarray}
where the source term is chosen as a linear function of the field variable
\begin{eqnarray*}
&& f(u) = \beta u(x,t) \ .
\end{eqnarray*}
The model describes the one-dimensional transport problem of a concentration $u(x,t)$
of a chemical or biological species in a flowing medium such as air or water.
The species concentration is assumed horizontally and vertically well mixed such that it varies only in the longitudinal or downstream direction. Moreover, a steady and uniform flow field is imposed and the effects of the dispersion are constant in time and space.
A reaction where the transformation rate $\beta$ is proportional to the species concentration is considered; according to the sign of rate, we may have either decay effects or not. 

\begin{figure}[!hbt]
\centerline{
\psfrag{x}{{$x$}}
\psfrag{t}{{$t$}}
\psfrag{U_j^n}{{$U_j^n$}}
\includegraphics[width=.7\textwidth]{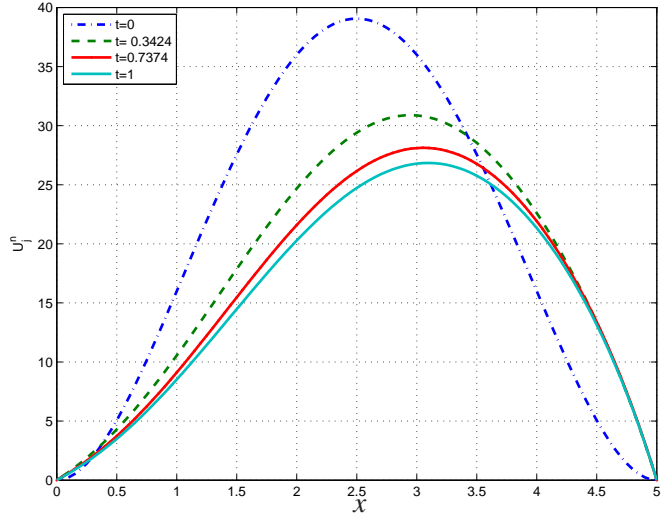}}
\caption{Numerical solutions of the FPDE (\ref{test2}) for $\alpha=0.5$ with $N=J=100$ at different time steps.}
\label{fig2}
\end{figure}

\begin{figure}[!hbt]
\centerline{
\psfrag{x}{{$x$}}
\psfrag{t}{{$t$}}
\psfrag{U_j^n}{{$U_j^n$}}
\includegraphics[width=.5\textwidth]{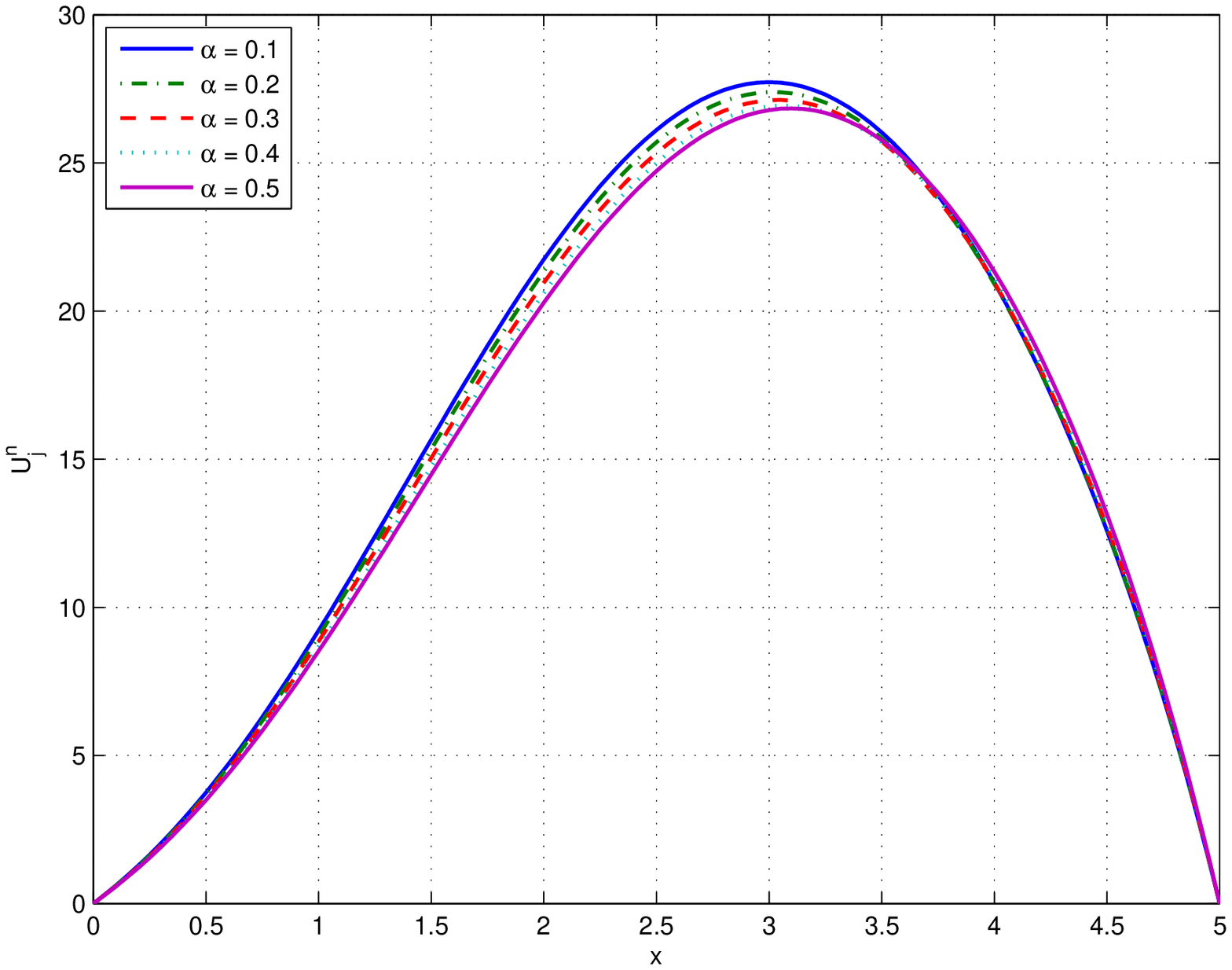}
\includegraphics[width=.5\textwidth]{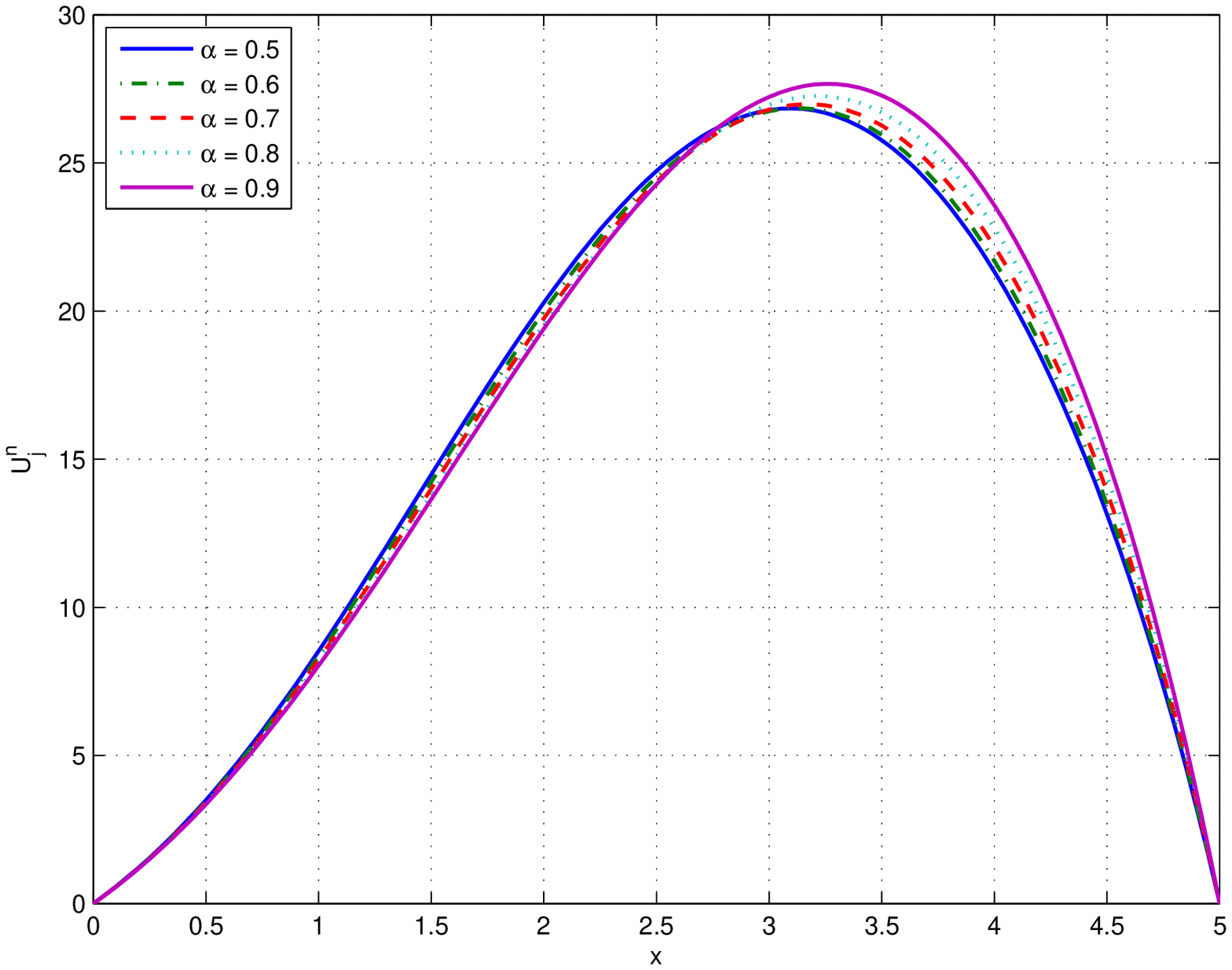}}
\caption{Numerical solutions of the FPDE (\ref{test2}) for different values of $\alpha$ with $N=J=100$.}
\label{fig2_2}
\end{figure}

In this example, we take $a=0$, $b=5$, $K_1=1$, $K_2=1$ and $\beta = 0.2$. Figure \ref{fig2} shows the solution
behavior at different times obtained for $\alpha = 0.5$ and with $N=J=100$. The height of the solution profile decreases as
the time increases.
Figures \ref{fig2_2} show the solution behavior with different values of $\alpha$ between $0$ and $1$ at the final time $t^N=1$.
Figures \ref{fig2_2} also show that the solution exhibits an anomalous diffusion behavior.
The height of the solution profile decreases when $0.1 \leq \alpha \leq 0.5$ and increases when $0.6 \leq \alpha \leq 0.9$.
The solution continuously depends on the $\alpha$-order of the time-fractional derivative.

\section{Concluding remarks}
In this paper, we develop an unconditionally stable finite difference method on quasi-uniform grid for solving
time-fractional advection-diffusion equations involving the Caputo fractional derivative.
The Caputo time derivative is discretized by means of a direct generalization of the well-known fractional
L1 formula (\ref{form_c}) to the case of quasi-uniform meshes.
The L1 formula takes into account the non-local character of the time-fractional operator and allows to improve the order 
of accuracy in time of the proposed method by using of the quasi-uniform time discretization obtained by the mesh (\ref{dt}).
We prove the stability and convergence of the proposed method. Numerical experiments are carried out to support the theoretical results.
The reported numerical experiments point out that the difference method is more accurate on the quasi-uniform grid than on the uniform mesh
and the convergence order is greater on the quasi-uniform mesh than on the uniform mesh.
Moreover, it is important to note that, in view of its simplicity, the method is applicable to a wide class of FADEs occurring in applied sciences.

\section*{Acknowledgements}
The research of this work was supported, in part, by the University of Messina and by the GNCS of INDAM.

\bibliographystyle{plain}


\end{document}